\documentclass[11pt]{article}
\usepackage[pctex32]{graphics}
\usepackage{amssymb}
\usepackage{latexsym}
\usepackage{epsfig} 
\usepackage{multirow} 
\usepackage{pstricks}
\usepackage{float}
\usepackage{epstopdf}
\usepackage{mathrsfs}
\usepackage{soul}

\usepackage{amsmath}
\usepackage{amssymb}
\usepackage{relsize}
\usepackage{multirow}
\usepackage{array}
\usepackage{arydshln}

\usepackage{graphicx}
\usepackage{multicol} 

\usepackage{times}
\usepackage{epstopdf}

\definecolor{Red}{rgb}{1,0.,0.}
\usepackage{fancyhdr}
\usepackage{setspace}
\usepackage{algorithm}
\newcommand{\R}{{\mathbb R}}

\newcommand{\mW}{{\mathsf W}}
\newcommand{\mH}{{\mathsf H}}

\newcommand{\mT}{{\mathsf T}}
\newcommand{\mC}{{\mathsf C}}
\newcommand{\mI}{{\mathsf I}}
\newcommand{\mA}{{\mathsf A}}

\newcommand{\mE}{{\mathsf E}}

\newcommand{\mD}{{\mathsf D}}

\newcommand{\mGamma}{{\mathsf \Gamma}}

\newcommand{\COtwo}{{\rm CO}_2}
\newcommand{\proton}{{\rm H}^+}
\newcommand{\carbacid}{{\rm H}_2{\rm CO}_3}
\newcommand{\bicarbonate}{{\rm HCO}_3^-}
\title{Bayesian dictionary learning estimation of cell membrane permeability from surface pH data }
\author{{A Bocchinfuso}\footnote{Current address: CINECA, consorzio interuniversitario} \and D Calvetti   \and E Somersalo}
\date{Department of Mathematics, Applied Mathematics and Statistics \\
Case Western Reserve University
}

\begin{document}
\maketitle

\begin{abstract}
Gas transport across cell membrane is a very important process in cell biology, as it is essential for many crucial tasks, including cell respiration and pH regulation. Since it is difficult to observe directly how gasses cross cell membranes, mathematical models have played and continue to play an important role in providing a context for the experimental findings depending on gas exchange, and to test different hypotheses about the gas passage mechanism. The first family of gas transport models has been able to match qualitatively the experimental data, and a second generation of computational models has reproduced and explained the dynamic range of the data. These models depend on parameters that cannot be measured experimentally, but can be estimated from the available data and some a priori knowledge about the problem based on first principles of cell physiology.  The membrane permeability is arguably one of the most important parameters, and can also be used also to shed light on the validity of different models of the transport mechanisms. In this contribution we propose a new, computationally efficient approach to the cell permeability estimation from measurements of the pH on the cell membrane over the course of a lab experiment. The methodology is inspired by a technique used in data science applications known as dictionary learning. \\
{\bf Relevance to life sciences. }
In the late 1990's, the suggestion that gasses are transported via preferred gas channels embedded into the cell membrane challenged the century-old Overton's theory that gases pass through the lipid cell membrane by diffusing across it driven by the concentration gradient. Since the laboratory experiments alone do not provide enough evidence to favor one or another of the proposed mechanisms, mathematical models have been introduced to provide a context for the interpretation of the measurements. In this article we address the problem of estimating intrinsic properties of the cell membrane that could help resolve the question about the role of aquaporins and other membrane-embedded proteins in gas transport across cell membrane from typical pH measurements on the membrane of the oocytes of Xenopus laevis.   \\ {\bf Mathematical content: } 
Following up on a previous work where the membrane permeability was estimated using particle filter techniques, in this article we propose an algorithm based on dictionary learning for estimating cell membrane permeability and other pertinent model parameters. Dictionary learning methods have proven to provide versatile tools to solve inverse problems in which the insufficiency of the available data or the complexity of the forward model makes the use of optimization-based inverse solvers impractical. 
Computed examples illustrate that the novel approach, which  can be applied when the properties of the membrane do not change in the course of the data collection process, is computationally more efficient than the particle filter approach.
\end{abstract}

\section{Introduction}

\begin{figure}
	\centering{\includegraphics[width=7cm]{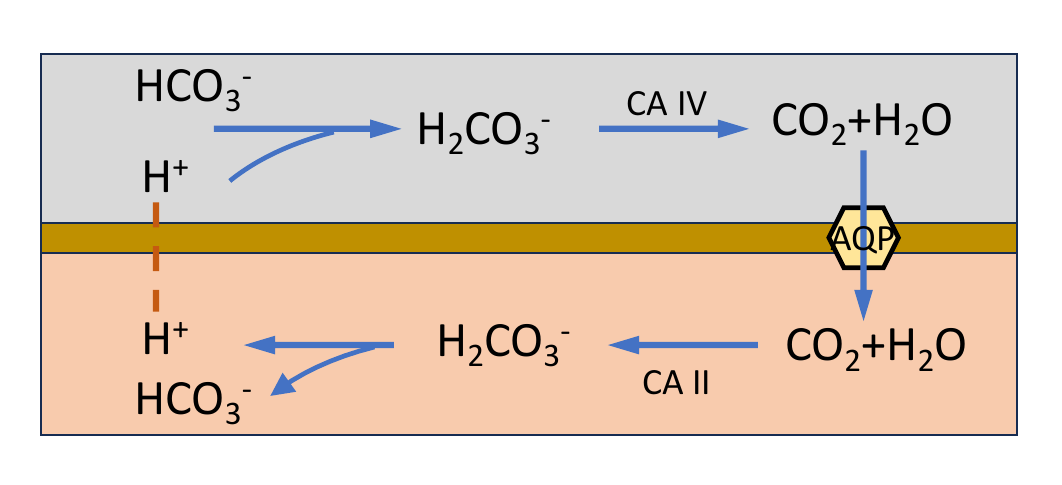}}
	\caption{{A schematic illustration of the $\COtwo$ based pH regulation across the cell membrane. A passage of a $\COtwo$ molecule, together with the reversible association/dissociation of the carbonic acid leads to a virtual transport of a proton across the membrane.
    The carbonic anhydrase (CA) enzyme acting on the membrane and inside the cell promote the dissociation or the association process of $\carbacid$.}}
	\label{fig:CA_action} %label{fig:shuttle}
\end{figure}

Gas transport across the cell membrane has been studied extensively because of its essential role in many cell functions, including cell respiration and pH regulation. In the literature on acid-base physiology, the regulation of the exchange of carbon dioxide ($\COtwo$) between the cell and its surrounding environment has been a recurrent topic, because the shuttle mechanism effectuating the passage, shown schematically in Figure~{\ref{fig:CA_action}}, is one of the ways for the cell to regulate the pH near the membrane surface. 

The mechanism of gas exchange through the cell membrane and the properties of membrane permeability have been investigated since the 1890's, starting with Overton's experiments and his theory \cite{Overton}, based on the hypothesis that the gases diffuse through the lipidic cell membrane driven by the concentration gradient across the membrane. In the late 1990s a new hypothesis emerged, suggesting with mounting experimental evidence that, in addition to diffusion, specialized gas channels, namely aquaporins (AQPs) and Rhesus proteins (Rh), might play a significant role in the exchange of gasses between cytoplasm and the outside environment \cite{Boron2010,Michenkova,Nakhoul,Endeward}. These gas channels, which are protein complexes embedded in membranes of various kinds of cells, were postulated to allow selectively the exchange of certain gases across the cell membrane \cite{MusaAziz2009,Geyer2013relativeAQPs,Geyer2013relativeRh}.

To validate the gas channel hypothesis, several experimental studies have been published on the passage of $\COtwo$ across the membrane of the oocytes of \emph{Xenopus laevis} \cite{Cooper,Nakhoul,MusaAziz2009}, a cell large enough to make it possible to inject it with heterologous RNA so as to have it express different kinds of AQPs and Rh proteins on its membrane, thus allowing to investigate the effects of the AQPs (or Rh) on the cell membrane permeability to gases \cite{Endeward}. 
The mathematical model proposed in \cite{Somersalo2012} describes the standard experiment with the oocytes of \emph{X. laevis}, in which the cell is immersed in an aqueous solution with high concentration of $\COtwo$, and the pH value directly above the surface is measured by means of a micro-sensor. The mathematical model, which  includes a biochemical description of the shuttle mechanism relating the passage of $\COtwo$ and changes in pH, provides a context for the interpretation of the experimental findings. In this model, it is assumed that the only chemical species that passes through the {\em X. laevis} cell membrane is $\COtwo$, while all other chemical species that are present in the solution cannot cross the membrane as it is de facto impermeable to them. This hypothesis is supported by earlier studies, see, e.g., 
\cite{Fei,Costa}.

The standard experiment comprises the measurement of the surface pH on the cell membrane by a pH-sensitive electrode placed against the membrane. At the beginning of the experiment, the data show a steep rise of the surface pH caused by the association of protons and bicarbonate into carbonic acid that breaks into water and $\COtwo$ replacing the carbon dioxide passing through the membrane into the cell. As the experiment proceeds, the $\COtwo$ transport slows down as the concentration inside the cell rises towards that of the outside, and this, together with the diffusion from the exterior domain drives the system to a new equilibrium.
The virtual proton transport into the cell is enhanced by the presence of carbonic anhydrases (CA), a family of enzymes that increase the dissociation of the $\carbacid$ right outside the cell membrane, while promoting the reverse reaction right under the membrane surface \cite{Boron2010,Michenkova,Occhipinti,Musa2014}, as illustrated schematically in Figure~\ref{fig:CA_action}. The concentration and spatial distribution of CA vary for different kind of cells, and it is known to be significantly high, for example, in  erythrocytes. 

\begin{figure}	\centering{\includegraphics[width=14cm]{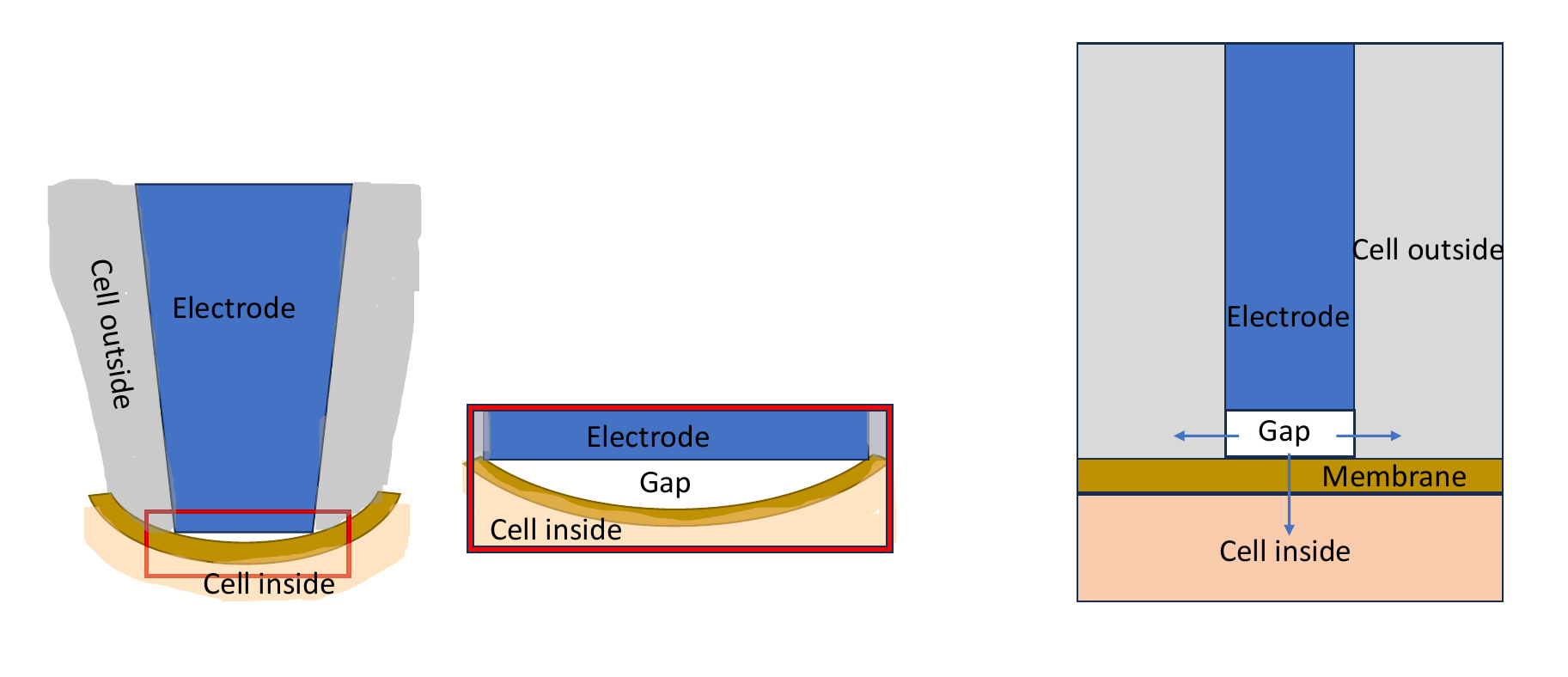}}
	\caption{{A schematic explanation of the hypothesized microenvironment created by the measurement device pushed against the membrane (left and center). The simplified geometry simulating the phenomenon is shown on the right. The free diffusion between the outside domain and the micro-compartment under the electrode is limited under the electrode edge.}}
	\label{fig:sensor_domain}
\end{figure}

The effect of carbonic anhydrases is part of the mathematical model proposed in \cite{Somersalo2012}, which reproduces well the qualitative behavior of pH during the experiment. However, the relatively simple spherically symmetric model fails to explain the observed dynamic range of the experimental data, the predicted pH peak values falling short of measured ones reported in the literature \cite{MusaAziz2009,Musa2014}. A plausible explanation of this shortcoming of the spherical model is that when pushed against the cell membrane, the measurement device itself creates a micro-environment, a nearly isolated pocket with biochemical dynamics different from the free membrane surface, see Figure~\ref{fig:sensor_domain}. To test this hypothesis, a detailed non-spherical finite element model for the microenvironment was developed in \cite{sensor}, and it was demonstrated that assuming a leaky clamping of the membrane patch under the electrode, was able to reproduce realistic dynamics in quantitative agreement with experiments.

To localize the effect of the membrane-bound carbonic anhydrase in the FEM model, it was necessary to resort to a dense spatial discretization, leading to a model with tens of thousands of degrees of freedom. Furthermore,  because of the very fast dissociation of carbonic acid and the slow diffusion process, the characteristic time scales range over several orders of magnitude, making the semi-discretized dynamic system extremely stiff. Consequently, to capture the fast-slow dynamics, in \cite{sensor} it was necessary to use an implicit-explicit scheme that separates different time scales. On a standard laptop or desktop computer, the resulting forward solver took typically several hours for each forward solve, making the algorithm impractical for estimating the membrane permeability and other parameters from the pH data. To overcome the issues of computational complexity, a reduced approximate model still accounting for the effect of the sensor was proposed in \cite{Modello}, and subsequently used for the estimation of the membrane permeability as the evolution step of a Bayesian particle filter algorithm in \cite{Filtro}. 

In recent years, there has been a growing interest in using dictionary learning methods for solving inverse problems when  the information in the data is not sufficient for the identification of any meaningful parametric model, or the forward model is not amenable for classical optimization-based approaches. Following this idea, here we propose a new alternative algorithm for solving the inverse problems of estimating the cell membrane properties based on a dictionary learning algorithm introduced in \cite{Dictionary}. The main advantage of the new approach is that it is computationally much more efficient than the particle filter algorithm, estimating the desired parameters to a similar accuracy as the particle filter at a small fraction of the computing time. 

In the following section, we review the details of the forward model that will be used to generate the dictionary model, sections 3 and 4 review the Bayesian dictionary learning algorithm.
The viability of the approach is demonstrated with computed examples in section 5.

\section{Experimental framework and mathematical model}

In this section we introduce the experimental setting used to estimate the cell membrane permeability to $\COtwo$ along with a number of other model parameters, and provide the corresponding mathematical description. 

The experiment starts with  a homogeneous solution
containing carbon dioxide ($\COtwo$), bicarbonate ($\bicarbonate$), carbonic acid ($\carbacid$), hydrogen ion ($\proton$), a non-specified pH buffer other than bicarbonate (${\rm A}^-$), and the relative acid (${\rm HA}$) in equilibrium with known initial concentrations. In the standard experiment, the {\em X. laevis} oocyte is immersed in this solution, and the micro-electrode measuring the membrane surface pH is already in place. The interior of the cell is assumed to consist of a solution of the same species as the exterior, but with different known equilibrium concentrations.

The chemical reactions tracked during the experiment are 
\begin{eqnarray}
	&& {\rm CO}_2 + {\rm H}_2{\rm O} {\mathop{\rightleftharpoons}^{k_1}_{k_{-1}} } {\rm H}_2{\rm CO}_3,  \label{react1}\\
	&&{\rm H}_2{\rm CO}_3 {\mathop{\rightleftharpoons}^{k_2}_{k_{-2}} } {\rm HCO}_3^- + {\rm H}^+, \label{react2}\\
	&& {\rm HA} {\mathop{\rightleftharpoons}^{k_{3}}_{k_{-3}} } {\rm A}^- + {\rm H}^+.\label{react3}
\end{eqnarray}
which occur both inside and outside the cell; the reaction rates $k_1$ and $k_{-1}$ of the first reaction are affected by the presence of carbonic anhydrases (CA) which have different concentrations and spatial distributions as explained later. 

As pointed out in \cite{sensor,Modello}, it is hypothesized that as the tip of the micro-electrode measuring the surface pH pushes against the membrane surface, a micro-environment under the tip different from the free membrane environment is created. Mathematically, the environment under the sensor tip is modeled as a pill-box cylindrical domain, as shown schematically in Figure~{\ref{fig:sensor_domain}}. While the rates of the reactions are not affected by the micro-environment, nor is the diffusion inside the domain, the
mobility of various species passing between the pill-box domain and the free membrane environment is reduced.

In \cite{sensor}, where the micro-environment under the tip of the sensor was modeled using a 3D finite element method approach, a substantial computational effort was required for each simulation, highlighting the need for a lower complexity model for the solution of the inverse problem of membrane permeability. Following \cite{Modello}, here we resort to a lumped compartment model for the perturbation caused by the sensor tip near the cell membrane. Based on the reasonable assumptions that the diffusion is mainly in the radial direction, and that the electrode has a minimal effect on the diffusion away from the membrane, the lumped model is effectively a radial diffusion model, ignoring for the most part diffusion in the direction parallel to the cell membrane, except right at the cell membrane, where diffusion from the free space to the micro-environment under the electrode tip is accounted for. The model is discussed in detail in the following paragraphs.

In our computational scheme, the reference frame is placed at the center of the oocyte, which is modeled as a sphere radius $R$. The changes in the concentration of the species are governed by a system of reaction-diffusion equations of the form 
\begin{equation}\label{3D}
	\frac{\partial u_n^\pm}{\partial t} = \nabla\cdot\kappa^\pm_n\nabla u_n^\pm + \sum_{j=1}^M S_{nj}\varphi^\pm_{j},
\end{equation}
where $u_n^\pm = u_n^\pm(x,t)$ is the concentration of the $n$th species at time $t$ and point $x$, the superscript $\pm$ indicates whether the concentration refers to that inside (-) or the outside (+) the cell membrane, $\varphi^\pm_{nj} = \varphi^\pm_{nj}(x,t)$  is the reaction flux of the $j$th reaction and, finally, $S_{nj}$ are the entries of the stoichiometric matrix $S$.  More specifically, the numbering of the concentrations is $u_1 = [{\rm CO}_2]$, $u_2 = [{\rm H_2CO}_3]$, 
$u_3 = [{\rm HCO}_3^-]$, $u_4 = [{\rm H}^+]$, $u_5 = [{\rm HA}]$, and $u_6 = [{\rm A}^-]$. The diffusion coefficients $\kappa_n^\pm$, assumed to be known, are set to the values reported in the literature listed in Table~\ref{tab:diffusion}.

\begin{table} [t]
	\begin{center}
		\begin{tabular}{l | l | c}
			\hline
			$R$ & Oocyte radius & $650\,\mu{\rm m}$ \\
			$R_\infty$ & External radius & $800\,\mu{\rm m}$  \\
			$w$ & Radius of the electrode tip & $10\,\mu{\rm m}$  \\
			$h$ & Distance of the electrode tip & $10\,\mu{\rm m}$  \\
		\end{tabular}
		\caption{\label{tab:geom} Geometric parameters of the cell and the electrode in the system. The parameter $R_\infty$ is the radius of the truncated computational domain. At the artificial boundary, Dirichlet boundary conditions corresponding to the initial equilibrium concentrations are imposed on the system.}
	\end{center}
\end{table}

\begin{table}[p]
	\begin{center}
	\begin{tabular}{l | c | c}
		Substance & Inside [$(\mu\,{\rm m})^2/s$] & Outside [$(\mu\,{\rm m})^2/s$] \\
		\hline
		${\rm CO_2}$ &  1.71$\times 10^3$  & 1.71$\times 10^3$ \\
		${\rm H_2CO_3}$ & 1.11$\times 10^3$ &1.11$\times 10^3$  \\
		${\rm HCO_3}^-$ & 1.11$\times 10^3$ &1.11$\times 10^3$  \\
		${\rm H}^+$  &  8.69$\times 10^3$  & 8.69$\times 10^3$ \\
		${\rm HA}$ &  1.56$\times 10^3$  & 1.56$\times 10^3$ \\
		${\rm A}^-$ &  1.56$\times 10^3$  & 1.56$\times 10^3$ \\
	\end{tabular}
	\caption{\label{tab:diffusion} Diffusion coefficients $\kappa$ of the substances in the solution inside and outside the cell membrane.}
	\end{center}

	\begin{center}
		\begin{tabular}{l | c | c}
			Substance & Inside [mM] & Outside [mM] \\
			\hline
			${\rm CO_2}$ &  0 & 0.4720 \\
			${\rm H_2CO_3}$ & 0 & 0.0013  \\
			${\rm HCO_3}^-$ & 0 & 9.901 \\
			${\rm H}^+$ (pH) & 6.310$\times 10^{-5}$ (7.2) & 3.162$\times 10^{-5}$ (7.5) \\
			${\rm HA}$ & 12.09 & 2.500\\
			${\rm A}^-$ & 15.22 & 2.500
		\end{tabular}
		\caption{\label{tab:initial} Initial concentration of the substances in the standard experiment. At the beginning of the simulation the oocyte is placed in the solution with the electrode in place. The distance between the cell and the electrode doesn't change through the simulation.}
	\end{center}
\end{table}

At the beginning of the experiment, the systems inside and outside the membrane are in equilibrium with concentrations constants throughout the exterior domain and within the cell, however, the constant concentrations are not equal across the membrane boundary.  The initial values of the concentrations are listed in Table~\ref{tab:initial}. The assumption that the solution is well-mixed and the concentrations of the species far from the cell are well-known justifies assigning Dirichlet boundary conditions at an artificial outer boundary of the spherical computational domain with radius $R_\infty$, 
\begin{equation}
	u_i(R_\infty,t) = u_i(R_\infty,0) = u_{i}^{\infty} \, , \; i=1, \dots, 6.
\end{equation}

The reaction fluxes both inside and outside the membrane are modeled according to the law of mass action,
\begin{eqnarray}
	\varphi_1(x,t)  &=&  k_1 u_1(x,t), \\ \label{mass action 1}
	\varphi_2(x,t)  &=&  k_{-1} u_2(x,t), \\ \label{mass action 2}
	\varphi_3(x,t)  &=&  k_2 u_2(x,t), \\ \label{mass action 3}
	\varphi_4(x,t)  &=&  k_{-2} u_3(x,t) u_4(x,t), \\ \label{mass action 4}
	\varphi_5(x,t)  &=&  k_3 u_5(x,t), \\ \label{mass action 5}
	\varphi_6(x,t)  &=&  k_{-3} u_6(x,t) u_4(x,t),\label{mass action 6}  
\end{eqnarray}
with the superscripts $\pm$ omitted to simplify the notation. The reaction rates valued used in our simulations are from the literature, see \cite{Modello} for details, listed in Table~\ref{tab:reaction rates}.

To account for the effect of CA, the fluxes $\varphi^\pm_1(x,t)$  and $\varphi^\pm_2(x,t)$ are enhanced by a multiplicative factor $A^\pm>1$ inside the cell and outside in the immediate surroundings of the cell membrane, yielding  
\[
\varphi_1^-(x,t) = A^- k_1 u_1^-(x,t), \quad \varphi_2^-(x,t) = A^- k_{-1} u_2^-(x,t), \quad   |x| < R
\]  
inside the cell, and 
\begin{equation*}
	\varphi_1^+(x,t)  =  \begin{cases} 
		A^+ k_1 u_1^+(x,t), \quad &\text{for } R<|x| <R+\delta, \\
		k_1 u_1^+(x,t), &\text{for } |x| \geq R+\delta,
	\end{cases} \quad 
	\varphi_2^+(x,t)  =  \begin{cases} 
		A^+ k_{-1} u_2^+(x,t), \quad &\text{for } R<|x| <R+\delta, \\
		k_{-1} u_1^+(x,t), & \text{for } |x| \geq R+\delta,
	\end{cases}
\end{equation*}
where $\delta$ is a small positive parameter defining the thickness of the
layer on top of the membrane enriched by the carbonic anhydrase, in agreement with \cite{Michenkova}.

\begin{table}
	\begin{center}
	\begin{tabular}{c | c | c | c}
		Reaction &  $k_\ell$  & $k_{-\ell}$ & $K=k_\ell/k_{-\ell}$ \\
		\hline
		${\displaystyle
			{\rm CO}_2 + {\rm H}_2{\rm O} {\mathop{\rightleftharpoons}^{k_1}_{k_{-1}} } {\rm H}_2{\rm CO}_3}$ &  0.0302  [1/s] & 10.9631  [1/s] & 2.7547$\times 10^{-4}$ \\
		${\displaystyle
			{\rm H}_2{\rm CO}_3 {\mathop{\rightleftharpoons}^{k_2}_{k_{-2}} } {\rm HCO}_3^- + {\rm H}^+}$  & $\varepsilon = 10^{-9}$ [1/s] & $\varepsilon/K_2$ &  $K_2 = 0.2407$ [mM]  \\
		$ {\displaystyle
			{\rm HA}_{\rm in} {\mathop{\rightleftharpoons}^{k_{3}}_{k_{-{3}}} } {\rm A}^-_{\rm in} + {\rm H}^+_{\rm in}}$ & $\varepsilon' = 10^{-6}$ [1/s]   &  $\varepsilon'/K_{\rm HA}$ &  $K_{\rm HA} = 7.9433\times 10^{-5}$  [mM]\\
		$ {\displaystyle
			{\rm HA}_{\rm out} {\mathop{\rightleftharpoons}^{k_{3}}_{k_{-{3}}} } {\rm A}^-_{\rm out} + {\rm H}^+_{\rm out}}$ & $\varepsilon' = 10^{-6}$ [1/s]   &  $\varepsilon'/K_{\rm HA}$ & $K_{\rm HA} = 3.1623\times 10^{-5}$ [mM]
	\end{tabular}
	\caption{\label{tab:reaction rates} Reaction rates. Note that $\varepsilon$ and $\varepsilon'$, whose precise values are not well known, represent the fast rates of the reactions for which no better guess is available. Furthermore, in the simulations we use $A k_{\pm 1}$, as those reactions are enhanced by the (multiplicative) acceleration factor $A$.}
\end{center}
\end{table}

In line with the assumption that $\COtwo$ is the only species crossing the cell membrane, the membrane permeability to all other species is set to zero, i.e., $\lambda_i = 0, \; i = 2, \ldots, 6$, so that the outside environment and the cell interior are coupled only through the $\COtwo$ flux. The rate at which $\COtwo$ passes the membrane is assumed to be proportional to the concentration gradient 
across the membrane, and the flux is modeled according to Fick's law,
\begin{equation}
	\kappa_1^+\frac{\partial u_1^+}{\partial r}(x,t)\bigg|_{r=R}  =  \kappa_1^-\frac{\partial u_1^-}{\partial r}(x,t)\bigg|_{r=R} = -\lambda \left(u^+_1(x,t)-u^-_1(x,t)\right)\big|_{r=R},
\end{equation}
where $\lambda$ is the permeability of the cell membrane to carbon dioxide. Observe that this model does not assume any particular mechanism for the gas transport across the membrane. The main hypothesis is that oocytes expressing the AQP and Rh proteins have notably higher permeability than the non-treated cells. Therefore, the inverse problem of estimating $\lambda$ from the pH data could provide important supporting evidence for the gas channel 
hypothesis \cite{MusaAziz2009,Musa2014}.

In line with the assumption that the  {\em X. laevis} oocyte is perfectly spherical with known radius, we write \eqref{3D} in spherical coordinates as
\begin{equation}
	\frac{\partial u_n^\pm}{\partial t} = \frac1{r^2}\frac{\partial}{\partial r}\left(\kappa^\pm_n r^2 \frac{\partial  u_n^\pm}{\partial r}\right) + \sum_{j=1}^M S_{nj}\varphi^\pm_{j}.
\end{equation}
The complete list of the parameters that specify the geometry of the system, including the radius of the cell and the radius of the sensor tip can be found in Table~\ref{tab:geom}.

The model described above is discretized in the spatial direction by one-dimensional first order finite elements, using a mesh with higher node density near the membrane to capture the local high gradients and the effect of CA, and the system of ordinary differential equations in time arising from the semidiscretization is solved using built-in stiff integrators of MATLAB. 

Following the ideas in \cite{Modello},
to account for the micro-environment under the electrode tip, a separate well-mixed compartment under the electrode is added to the model. The
concentration distributions in this small compartment are lumped into 
spatially averaged time-dependent concentrations $u_n^0(t)$. The compartment is coupled with the interior of the cell through $\COtwo$ transport across the membrane, as well as with the exterior domain through restricted diffusion through the partial clamping by the electrode rim. The dynamics of the concentrations in the domain is given by the equations 
\begin{equation}\label{eq:sensor}
	\frac{d  u_n^0}{d t}(t) = \delta_{n,1} \frac{\lambda}{h} \left(u^0_n(t) - u^-_n(R,t)\right) - \frac{2\gamma}{P}
	\left( u_n^+(R,t) -  u_n^0(t)\right)
	+ \sum_{j=1}^m S_{mj} \varphi_j^0(t).
\end{equation}
where $\lambda$ is the permeability of the cell to $\COtwo$, $h$ is the distance between the sensor tip and the membrane, and the parameter $\gamma$, referred to as the \emph{quenching factor}, accounts for the fact that the partial clamping of the electrode on the membrane reduces the mobility of the species between the micro-environment created by this interaction and the free environment. Finally, $P$ is the radius of the sensor tip. Here, we use the Kronecker symbol $\delta_{n,1}$ to indicate that the membrane is impermeable to the species $n\neq 1$.

The idealized data consist of measurements of the pH in the micro-environment,
\begin{equation}
	y_i = \log_{10} u^0_4 (t_i),  
\end{equation}
with sampling frequency of 2 Hz, i.e., $t_i = i\Delta t$ with $\Delta t = 0.5\,{\rm s}$, and $0\leq i\leq 1000$, that is, the duration of the experiment is 500 seconds. To incorporate an uncertainty to the simulated measurement, we assume that the device has a precision $p = 0.02$ pH units, and the pH reading is modeled as a quantized round-off value defined as 
\begin{equation}\label{quantization}
 y_i^{(q)} = \left[\frac{y_i}{p}\right] p,
\end{equation}
i.e., we round the simulated pH value to the nearest multiple of the precision $p$. Following \cite{Filtro}, the model parameters of interest to be estimated from the data are the membrane permeability $\lambda$, the CA enhancement factor $A^0$ at the outer surface of the membrane in the micro-environment, and the quenching parameter $\gamma$ controlling the mobility between the exterior domain and the micro-environment.
To set up the inverse problem, we non-dimensionalize the model parameters of interest by setting
\begin{eqnarray}
	\lambda &=& \xi_\lambda \, \lambda_0, \quad 0\leq \xi_\lambda \leq 1 \, , \\
	A^0 &=& \xi_A\, A, \quad    0\leq \xi_A \leq 1 \, , \\
	\gamma &=& \gamma_0 \, 10^{[-0.5 \, \xi_{\gamma} + (1 - \xi_{\gamma}) ]}, \quad 0\leq \xi_\gamma\leq 1.
\end{eqnarray} 
where $\lambda_0= 34.2\,\mu{\rm m}/{\rm s}$ corresponds to the reference value obtained by assuming that the membrane is non-resistant,  $A=20$ is an assumed value on the free surface, and $\gamma_0 = 10^{-4}\,\mu{\rm m}/{\rm s}$.
The inverse problem considered here is to estimate the dimensionless vector $\xi = (\xi_\lambda,\xi_A,\xi_\gamma)$ from the time series data $\big(y_i^{(q)}\big)_{i=0}^{1000}$.

\section{Bayesian dictionary learning}
\label{sec:dictionary}

In data science, dictionary learning and dictionary matching refer to a class of algorithms aiming at explaining a given datum in terms of a collection of precomputed or previously measured and labeled entries, referred to as the underlying dictionary. In applications where the new datum is a vector $b$ of  measured response of a physical system, dictionary learning can be employed to explain it as a linear combination of a collection of possible responses. Assuming that the known responses are stored as vectors $d_1, \dots, d_p$ of the same dimension as $b$, often referred to as atoms, they can be organized as the columns of a matrix dictionary $\mD$, and the problem can be formulated as finding a solution of the linear system 
\begin{equation}
	\mD x = b. \label{eq:linear_system}
\end{equation}
In inverse problem applications, the atoms may represent precomputed outputs corresponding to a model that depends on model parameters that one is primarily interested in, the model parameter values serving as labels. Identifying one or few atoms that are able to explain the datum allows an identification of the model parameters, thus providing a solution to the inverse problem.
The approach is particularly useful for problems in which the forward model is not amenable for classical methods, 
see, e.g., \cite{Ma,Badger} for examples. In the current application, the dictionary consists of precomputed pH curves, each corresponding to a given parameter value combination.

In order to identify model parameters corresponding to the observation, it is desirable to explain the new datum in terms of only few atoms, or equivalently, to have most of the entries of the solution vector $x$ vanish. This particular instance of dictionary learning is referred to as sparse coding.  Furthermore, 
in order for the datum to be interpreted as the superposition of atoms without cancellation of contributions from different atoms, the entries of the solution of (\ref{eq:linear_system}) should satisfy a non-negativity constraint $x \geq 0$, the inequality understood component-wise. Moreover, dictionary matching is often combined with dictionary reduction, or "learning the dictionary", aiming at a compressed dictionary that contains the salient features of the original one but leads to less demanding computations. 

For computational efficiency, it is often advantageous to preprocess the dictionary and identify possible clusters, i.e., divide the atoms in groups, each group consisting of atoms with similar characteristics. In lack of better criteria, one con resort to unsupervised clustering techniques such as $k$-means or $k$-medoids algorithm.
Assuming that the dictionary atoms are partitioned into $K$ clusters, refereed to as sub-dictionaries, following \cite{Dictionary}, we can permute the columns of $\mD$ so that columns corresponding to atoms in the same cluster are adjacent, yielding the block structure 
\begin{equation}
	\mD = \begin{bmatrix}
		\begin{array}{c:c:c}
			\mD^{(1)} &\dots &\mD^{(K)}
		\end{array}
	\end{bmatrix}. 
\end{equation}
In this manner the dictionary learning problem can be solved by first identifying the sub-dictionaries that may be pertinent, and subsequently reformulating (\ref{eq:linear_system}) with a smaller dictionary matrix that does not include the blocks corresponding to irrelevant sub-dictionaries. It has been shown in \cite{Dictionary} that this approach can substantially reduce the computational complexity, especially when the dictionary comprises a large number of atoms. 

In data science applications, it is customary to reduce the size of large data matrices by 
means of low-rank matrix approximations. If the columns of a data matrix have some degree of similarity, it is natural to identify a few vectors capturing their main characterizing features, and to approximate them as linear combinations of the feature vectors. This leads to a low rank representation of the data vectors.
Once the sub-dictionaries have been identified, e.g., by clustering, each one can be compressed and approximated by just a few feature vectors. More specifically, if $\mD^{(i)} \in\R^{m\times p_i}$, we write a low-rank approximation
\begin{equation}
	\mD^{(i)} \approx \mW^{(i)} \mH^{(i)},  \mbox{ where $\mW^{(i)} \in\R^{m\times k_i}$, $\mH^{(i)} \in\R^{k_i\times p_i}$, and  $k_i<p_i$,}
\end{equation}
that is, every atom in the sub-dictionary $\mD^{(i)}$ is approximated by a linear combination of the columns of $W^{(i)}$, referred to as the $i$th code book. The rank $k_i$ represents the effective dimensionality of the subdictionary, and  $W^{(i)}$ represents a reduced sub-dictionary because it summarizes $\mD^{(i)}$ in terms of much fewer feature vectors. Observe that clustering the the dictionary entries in subdictionaries of similar entries reduces the diversity within each group, thus helping to lower the number of feature vectors needed to explain each subdictionary.
A popular low-rank approximation technique, the principal component analysis (PCA), is based on the singular value decomposition (SVD) of the matrix, however, if positivity constraints for the feature vectors and/or the coefficients are necessary,  non-negative matrix factorization (NMF) algorithms may be preferable \cite{Gillis2020nonnegative}.
While for a given approximation rank $k_i$, the approximation properties of NMF measured in the Frobenius norm are inferior to those of PCA, the non-negativity of the feature vectors and of the weights facilitate the interpretation of the factorization.  

In the application of interest here, the atoms and the data are nonnegative vectors, since they consist of time series of pH measurements; furthermore, it is known \cite{Somersalo2012,Modello,Filtro,sensor} that in the standard experiment in which the $\COtwo$ concentration in the exterior domain is increased above the corresponding concentration inside the cell, the pH value never falls below the reference value ${\rm pH}_0$ corresponding to the initial equilibrium, which in the current simulation is ${\rm pH}_0 = 7.5$, temporarily rising above it and subsequently decreasing towards it as the system returns to equilibrium. In line with this observation, the equilibrium value is subtracted from the quantized measurements $y_j^{(q)}$, and the data used in our computations, defined as a vector
\begin{equation}\label{remove bg}
	b = y^{(q)} - {\rm pH}_0, 
\end{equation}
has all components nonnegative. 
To respect the non-negativity, for each subdictionary we compute a corresponding NMF of the form 
\begin{equation}\label{NMF}
	\mD^{(i)} \approx \mW^{(i)} \mH^{(i)}, \qquad
	\mW^{(i)}_{l,m} \geq 0, \;
	\mH^{(i)} _{l,m} \geq 0 \, ,
\end{equation}
the inequalities being understood component-wise.
The columns of the matrix $\mW^{(i)}$, sometimes referred to as feature vectors, can be interpreted as a summary of the atoms of $\mD^{(i)}$, thus $\mW^{(i)}$ can be thought of as a compressed version of $\mD^{(i)}$. In the current application, we compute the NMF using  an iterated nonnegative least squares algorithm with sparsity constraints for the matrix $\mH^{(i)}$ as discussed in more detail later.

In general, the approximation (\ref{NMF}) is inexact,
and to model the discrepancy, we write
\begin{equation}
	\mD^{(i)} = \mW^{(i)} \mH^{(i)} + \mE^{(i)}\, ,
\end{equation}
where the discrepancy matrix $\mE^{(i)}$ can be explicitly computed.
In the proposed algorithm, we adopt a probabilistic view, thus the columns of the matrix $\mD^{(i)}$, and therefore those of $\mE^{(i)}$, are modeled as independent samples from some underlying probability distribution. In turn, if the columns of $\mE^{(i)}$ are independent realizations of a random variable, denoted by $E^{(i)}$,
we can estimate the underlying probability distribution from a sample of realizations of this random variable.  Therefore, for each subdictionary, we  consider the sample of dictionary entries $\{d_j^{(i)}\}_{j=1}^{K_i}$ comprising the columns of the matrix $\mD^{(i)}$,
find their approximate representations in terms of the feature vectors, and compute the approximation errors,
\begin{equation}\label{error sample}
e^{(i)}_j = d^{(i)}_j - \mW^{(j)}h^{(i)}_j\,\quad 1\leq j\leq K_i\, .
\end{equation}

For computational convenience, the probability density of this random variable is approximated by a normal distribution,
\[
	E^{(i)} \sim \mathcal{N}(\mu^{(i)}_{\rm DCE}, C^{(i)}_{\rm DCE}),
\]
where the mean $\mu^{(i)}_{\rm DCE}$ and covariance $\mC^{(i)}_{\rm DCE}$ are estimated from columns of the matrix $\mE^{(i)}$, that is,
\begin{equation}\label{DCE mean}
 \mu_{\rm DCE}^{(i)} = \frac {1}{K_i}\sum_{j=1}^{K_i} e_j^{(i)},
\end{equation}
\begin{equation}\label{DCE cov}
\mC_{\rm DCE}^{(i)} = \frac 1{K_i-1}\sum_{j=1}^{K_i}\big(e_j^{(i)} - \mu_{\rm DCE}^{(i)}\big)
\big(e_j^{(i)} - \mu_{\rm DCE}^{(i)}\big)^\mT + \delta\,\mI_m\, ,
\end{equation}
where $\delta>0$ is a small regularization parameter to guarantee the positive definiteness of the estimated covariance matrix.
The details of the sample generation are discussed further in connection with the computed example.

Given a time series $b$ of pH measurements on the cell membrane during an experiment and a collection of subdictionaries $\mD^{(1)}, \ldots, \mD^{(K)}$, the first step in the estimation of the membrane permeability is to identify the subdictionary $\mD^{(\hat{i})}$ that best explains $b$. The best explaining dictionary is defined here as the one for which the solution of the whitened linear least squares problem with non-negativity constraint
\begin{equation}\label{local LSQ}
h^{(i)} = \underset{h\geq 0}{\rm argmin} \left\| \left(\mC^{(i)}_{\rm DCE}\right)^{-1/2} \left(b -\mu^{(i)}_{\rm DCE} - W^{(i)} h\right) \right\|
\end{equation}
yields the smallest residual, that is, the index $\widehat i$ of the best explaining dictionary is
\begin{equation}\label{winner}
\widehat{i} = \underset{i}{\rm argmin} \left\| \left(\mC_{\rm DCE}^{(i)}\right)^{-1/2} \left(b -\mu^{(i)}_{\rm DCE} - \mW^{(i)} h^{(i)}\right) \right\| \, .
\end{equation}

Once the best subdictionary has been identified, the datum $b$ is interpreted in terms of its atoms by solving the non-negative linear least squares problem 
\begin{equation}\label{eq:LS_reconstruction}
\mbox{
minimize		$\big\|\mD^{(\widehat{i})} x - b\big\|$ subject to $x\geq 0$,  $\| x \|_0 \ll p_{\,\widehat{i}}$,
}
\end{equation}
where $p_{\,\widehat{i}}$ is the number of atoms in $\mD^{(\widehat{i})}$, and $\|x\|_0$ refers to the cardinality of the support of the vector $x$. This minimization is done by means of the sparsity promoting algorithm described in subsection \ref{sub:IAS}. We denote the sparse solution of the problem (\ref{eq:LS_reconstruction}) by $\widehat x$.

The final step of the process entails going from the solution  vector $\widehat x$  to the parameter vector $\xi$ associated with the measured data $b$.
If $\widehat x$ contains only a single non-zero component, the parameter identification is straightforward. When the sparse vector $x$ contains more than one non-zero entry, a winner-takes-all strategy of electing the parameter vector corresponding to the best matching dictionary atom is a feasible approach.
In the computed examples of this paper, we use an approximate interpolation scheme.
Consider the mapping from the parameter space to the data space, $f:\xi\mapsto d$. We denote by $\xi_j^{(\hat i)}$ the parameter vectors corresponding to the best matching subdictionary, $1\leq j\leq p_{\,\hat i}$, and by $d_j^{(\hat i)}$ the corresponding computed outputs.
We normalize the coefficients $\widehat x_j$, defining
\[
\eta_j = \frac1{\|\hat x\|_1}\hat x_j,
\]
and define the convex combination of the parameter vectors,
\[
\overline\xi = \sum_{j=1}^{p_{\hat i}} \eta_j\xi_j^{(\hat i)}
\]
Assuming that the individual parameter vectors $\xi^{(\hat i)}_j$ are close to the convex combination $\overline \xi$, we obtain by linearization around $\overline\xi$ the approximate identity
\begin{eqnarray*}
b&\approx&\sum_{j=1}^{p_{\,\hat i}} \widehat x_j d^{(\hat i)}_j = \frac 1{\|\hat x\|_1} \sum_{j=1}^{p_{\,\widehat i}}\eta_j f\big(\xi^{(\hat i)}_j\big) \\ &\approx&  
\frac 1{\|\hat x\|_1} \sum_{j=1}^{p_{\,\widehat i}}\eta_j\left(f(\overline \xi) + Df(\overline \xi)(\xi^{(\hat i)}_j - \overline\xi)\right) \\
&=&\frac 1{\|\widehat x\|_1} \left(f(\overline \xi) +Df(\overline \xi)\big(\sum_{j=1}^{p_{\,\hat i}}\eta_j \xi^{(\hat i)}_j - \overline\xi\big)\right) \\
&=& \frac 1{\|\widehat x\|_1} f(\overline\xi).
\end{eqnarray*}
This calculation indicates that the best matching parameter value explaining the data is given by the convex combination
\begin{equation}\label{interpolate}
 \widehat \xi = \sum_{j=1}^{p_{\,\hat i}} \eta_j\xi_j^{(\hat i)}.
\end{equation}

\subsection{Bayesian sparsity promoting linear least squares solvers}
\label{sub:IAS}
The Bayesian dictionary learning procedure outlined above relies on the an algorithm for computing sparse and nonnegative solutions of linear least squares problems. In the following we provide a brief overview of the Iterative Alternating Sequential (IAS) algorithm, For the derivation of the IAS algorithm and a detailed discussion of its properties we refer to \cite{BSC}. 

Consider a linear inverse problem with non-negativity constraint of the form 
\begin{equation}\label{eq:LS_system}
	b = \mA x + \varepsilon, \quad x\geq 0, 
\end{equation}
where $\mA\in\R^{m\times p}$ is a matrix encoding the relationship between the unknown input $x$ and the measured data $b$ and $\varepsilon$ is an additive observation error modeled as a zero-mean Gaussian random variable with covariance $\mGamma$. Under these assumptions, the likelihood is of the form 
\[
	\pi_{B\mid X}(b\mid x) \propto \exp \left( -\frac{1}{2} (b - \mA x)^\top \mGamma^{-1} (b - \mA x) \right).
\]
Without a loss of generality, we may assume that the data are whitened, i.e., $\mGamma = \mI$, the identity matrix.
In the Bayesian framework the unknown $x$ is modeled as a random variable and we express the belief that it should have only a few nonzero entries by means of a zero-mean conditionally Gaussian prior of the form 
\[
	X \sim \mathcal{N}(0, \mD_\theta) \, , \;  \mD_\theta=\text{diag}(\theta_1, \dots, \theta_p) \, , \; \theta_j > 0,
\]
or,
\[
	\pi_{X\mid \Theta}(x\mid\theta) \propto \pi_+(x) {\rm exp}\left( -\frac{1}{2} \sum_{j=1}^{p} \log \theta_j - \frac{1}{2} \big\| \mD_\theta^{-1/2} x \big\|^2 \right),
\]
where $\pi_+(x)$ is the characteristic function of the non-negative cone.

The unknown variances $\theta_j$ are modeled as random variables. Furthermore, assume that they are mutually independent 
with a distributions that promote small values while allowing large outliers. For computational convenience we assume that the variances are distributed according to a gamma distribution, 
\begin{equation}
	\pi_{\Theta}(\theta) \propto \prod_{j=1}^{p} \frac{1}{\vartheta_j^{\beta}} \theta_j^{\beta-1} \exp\left(-\frac{\theta_j}{\vartheta_j}\right),
\end{equation}
where $\beta$ is the shape parameter and $\vartheta_j$ is the scale parameter.
For an argument supporting the sparsity assumption, see \cite{L2 Magic,Dictionary,BSC} for further discussion.

It follows from Bayes' Theorem and the law of total probability that the posterior distribution of the pair $(X,\Theta)$ is of the form 
\begin{equation}
	\pi_{X,\Theta\mid B}(x, \theta\mid b) \propto \pi_+(x)\exp \left( -\frac{1}{2} \Vert b - \mA x \Vert^2 - \frac{1}{2} \big\| \mD_\theta^{-1/2} x \big\|^2 + \bigg( \beta - \frac{3}{2} \bigg) \sum_{j=1}^{p} \log\theta_j - \sum_{j=1}^{p} \frac{\theta_j}{\vartheta_j} \right) \, . \label{eq:IAS}
\end{equation}

One way to summarize the posterior is to compute  the Maximum a Posteriori (MAP) estimate, which is the realization maximizing the posterior probability density. The block descent IAS algorithm for the MAP estimate computation proceeds by alternatively updating $x$ and $\theta$, organizing the computations as indicated in the following algorithm.

\bigskip
\hrule

{\bf IAS algorithm with non-negativity constraint}
\medskip

\hrule
\begin{itemize}
\item[] {\bf Given:} matrix $\mA\in\R^{m\times p}$, $b\in\R^m$, $\beta$, $\vartheta_j, \; 1 \leq j \leq p$,
\item[] {\bf Initialize:}  $\theta^{(0)}_j = \vartheta_j, \; 1 \leq j \leq p$, set the counter $k=0$.
\item[] {\bf Iterate} until stopping criterion is met:
\begin{enumerate}
\item Update $x$:  Let $\theta_j = \theta^{(k)}_j, \; 1\leq j \leq p$ and set 
\[
x^{(k+1)} = \underset{x\geq 0}{\rm argmin} \big\{ \| b - \mA x\|^2 + \big\|\mD_{\theta^{(k)}}^{-1/2} x \big\| ^2 \big\}.
\]
\item Update $\theta_j, \; 1 \leq j \leq p\;$by solving
\[
 \frac{\partial {\mathscr E}}{\partial\theta_j}(x^{(k+1)}, \theta) = 0,
\]
yielding 
\[
\theta_j^{(k+1)} =  \frac{\vartheta_j}{2}\left( \eta + \sqrt{\eta^2 + \frac{[x_j^{(k+1)}]^2}{4\vartheta_j}}\right), \quad \eta = \beta - \frac 32.
\]
\item Advance the counter $k\to k+1$ and check convergence.
\end{enumerate}

\bigskip
\hrule
\bigskip
\end{itemize}

The properties of the IAS algorithm, and its generalizations using generalized gamma distributions have been investigated and reported in the literature, see \cite{BSC} for a summary. In particular, without the non-negativity constraint, it has been shown that when the gamma distribution for the variances is used, the objective function minimized by the IAS algorithm is strictly convex, hence it has a unique global minimum. Under those conditions the IAS iterates converge to the unique minimizer and moreover, in the limit as $\eta = \beta - 3/2 \rightarrow 0^+$, the solution computed by he IAS algorithm converges to a weighted $\ell_1$-regularized solution, 
\[
 x_{\ell_1} = {\rm argmin} \left\{ \frac 12 \|b - \mA x\|^2 + \sqrt{2} \sum_{j=1}^p \frac{|x_j|}{\sqrt{\vartheta_j}}\right\},
\]
This limiting argument confirms that for small $\eta$, the algorithm favors sparse solutions. Moreover, the presence of the hyperparameters $\vartheta_j$ in the above formula relates them to the classical Tikhonov regularization parameters as discussed in detail in  \cite{calvetti2025distributed}.
Criteria for selecting the values of the hyperparameters $\vartheta$ by a sensitivity analysis argument can be found in  \cite{L2 Magic}. The non-negativity constraint for the IAS was introduced in \cite{Inverse_gamma}, where each updating step of $x$ by a least squares solution is followed by a projection step to the non-negative cone.

The IAS algorithm has been proved to be flexible and easy to modify to accommodate the needs of the problem in question, as it reduces the sparsity promoting inversion in two technically simple steps, a least squares problem with possible constraints for updating $x$, and a function evolution for updating $\theta$. The formulation of the first step as a standard least squares problem makes it possible to take advantage of the power of numerical linear algebra for the solution of linear least squares problems. Krylov subspace least squares solvers are the methods of choice for  very large problems, including those in which the matrix cannot be stored in the memory, or when the linear operator is defined in a matrix-free form.

In \cite{Waniorek}, the IAS algorithm was proposed for finding sparse NMF factorizations for non-negative matrices, by reducing the NMF algorithm in an alternating sequence of non-negative least squares problems, see \cite{Gillis2020nonnegative}.
When using the IAS algorithm for updating the low-rank factors of the nonnegative matrix factorization, it is necessary to ensure that the computed solution is nonnnegative. The nonegativity is enforced by solving the linear least squares problem yielding the updated $x$ using the CGLS method, projecting the computed solution onto the nonnegative cone, and using it as the initial guess for restating the CGLS method, see
\cite{calvetti2020sparsity}.
In our application, we restart the CGLS method only once. 

It has been observed that sparsity promotion is stronger when the $\theta$ components follow an inverse gamma distribution, although the non-convexity of the resulting objective function could cause the IAS algorithm to stop at local minima. To combine the global convergence of the gamma hyperprior with the stronger sparsity promotion of the inverse gamma hyperprior, we start the IAS algorithm with the gamma hyperprior, and after a certain number of iterations once we are close to the global minimizer, we switch to the inverse gamma. The computational details of this hybrid scheme can be found in \cite{Inverse_gamma}.

\section{Bayesian Dictionary Parameter Estimation Workflow}

In this section we summarize the different tasks that need be completed to design a dictionary-learning based parameter estimation algorithm for problems where the forward model is available and a bounded domain in the parameter space where the parameters to be estimated live is known.

\subsection{Dictionary preparation and reduction}

The first step is the creation of the dictionary, consisting of the output of simulations corresponding to different combinations of parameters. The discretization of parameter space, in addition to considerations of computational complexity,  should be chosen on the basis of the desired resolution of the parameter estimate. Denoting by $\{\xi_j\}_{j=1}^p$ the discrete set of parameter values covering the parameter domain and by $f$ the forward map, we generate simulations $\{d_j\}_{j=1}^p$, where $d_j = f(\xi_j)$. The simulations are assembled in the dictionary matrix $\mD =[d_1,\cdots,d_p]$. 

To reduce the computational cost of the procedure, the dictionary atoms are partitioned into $K$ clusters using either a partitioning of parameters space or by means of an unsupervised clustering algorithm. We reorganize the dictionary so that output vectors in the same class are adjacent to each other,
\[
\mD = \begin{bmatrix}
					\begin{array}{c:c:c}
						\mD^{(1)} &\dots & \mD^{(K)}
					\end{array}
				\end{bmatrix} ;
\]

Next we apply a low-rank factorization $ \mD^{(i)} = \mW^{(i)} \mH^{(i)} $ to each subdictionary. The rank $k_i$ of the approximation may vary from subdictionary to subdictionary.  

For each subdictionary, we compute the Dictionary Compression Error statistics in terms of the mean and covariance given by formulas (\ref{DCE mean}) and (\ref{DCE cov}).

{\bf Phase 1: Subdictionary Identification}

Given a data vector $b$, if the corresponding parameter vector were among those used to construct the dictionary, the search for the parameter estimate could be restricted to a single subdictionary. In Phase I, we therefore seek to express $b$ in terms of the feature vectors $\mW^{(i)}$, by solving the least squares problem (\ref{local LSQ}) for each subdictionary, and select the best matching subdictionary by the criterion (\ref{winner}).

{\bf Phase 2: Dictionary coding with deflation}

Once the best matching subdictionary $\mD^{(\widehat i)}$ has been identified in Phase 1, we proceed to coding $b$ in terms of the deflated dictionary consisting only of the subdictionary $\mD^{(\widehat i)}$ by solving the linear system 
\[
\mD^{(\widehat i)} x = b
\]
augmented by the sparsity constraint, in the sense (\ref{eq:LS_reconstruction}). The implementation of this step is done by using the IAS algorithm with positivity constraint. Denote by $\widehat x$ the sparse solution.

{\bf Phase 3: From dictionary coding to parameter space}

Having the sparse reconstruction of the data in terms of the entries of a single subdictionary, we interpolate the parameters to have the final estimate given by formula (\ref{interpolate}).  

\section{Computed  example}\label{sec:numerical}

In this section, the dictionary learning approach is applied to the inverse problem described in section 2. We start by sampling the parameter space in order to generate the dictionary: The feasible region for the unknown parameters,
\begin{equation}
	(\xi_\lambda, \xi_A, \xi_\gamma) \in [0.6, 1] \times [0.6, 1] \times [0, 1] \label{eq:grid}
\end{equation}
is sampled through a uniform grid of size $40 \times 40 \times 40$,
and for each node of the grid we compute the time series pH data corresponding to the measurement.  More specifically, the pH is computed for $t_j = j\delta t$, where $0\leq j\leq 1\,000$ and $\Delta t = 0.5\,{\rm s}$, using the electrode model, quantized with precision $p=0.02$ as in equation (\ref{quantization}), and the background removed as in (\ref{remove bg}). Each  vector containing the times series of pH values is an atom in the dictionary, yielding a dictionary of size $1\,001\times 64\,000$.

Following the proposed workflow, we partition the dictionary into $k$ classes using the $k$-medoids algorithm, where the distance matrix is based on the Euclidean distance. While the selection of the number of clusters in unsupervised clustering is somewhat arbitrary, we use the scree plot criterion (or "elbow criterion") for the selection, see, e.g., \cite{444}, suggesting the cluster number $k=7$.

\begin{figure}
    \centerline{
\includegraphics[width=9cm]{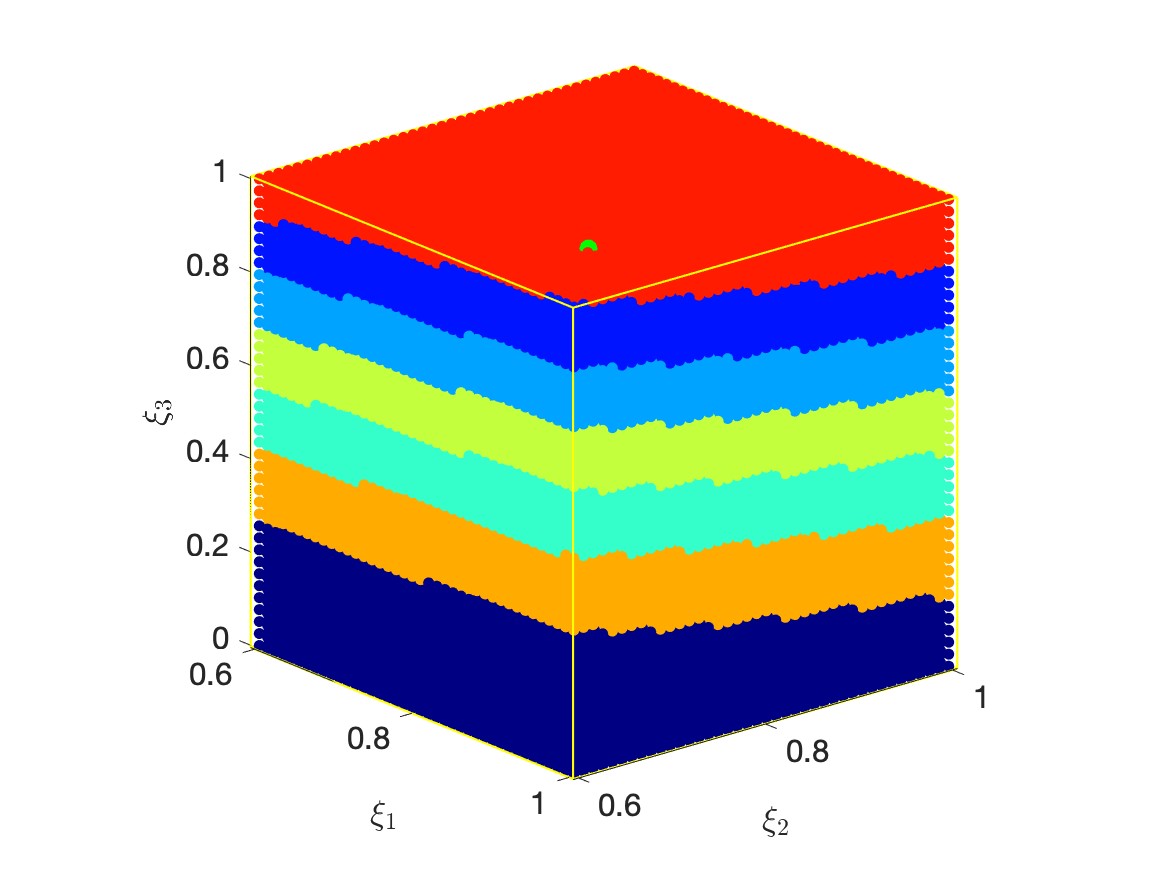}
 \includegraphics[width=9cm]{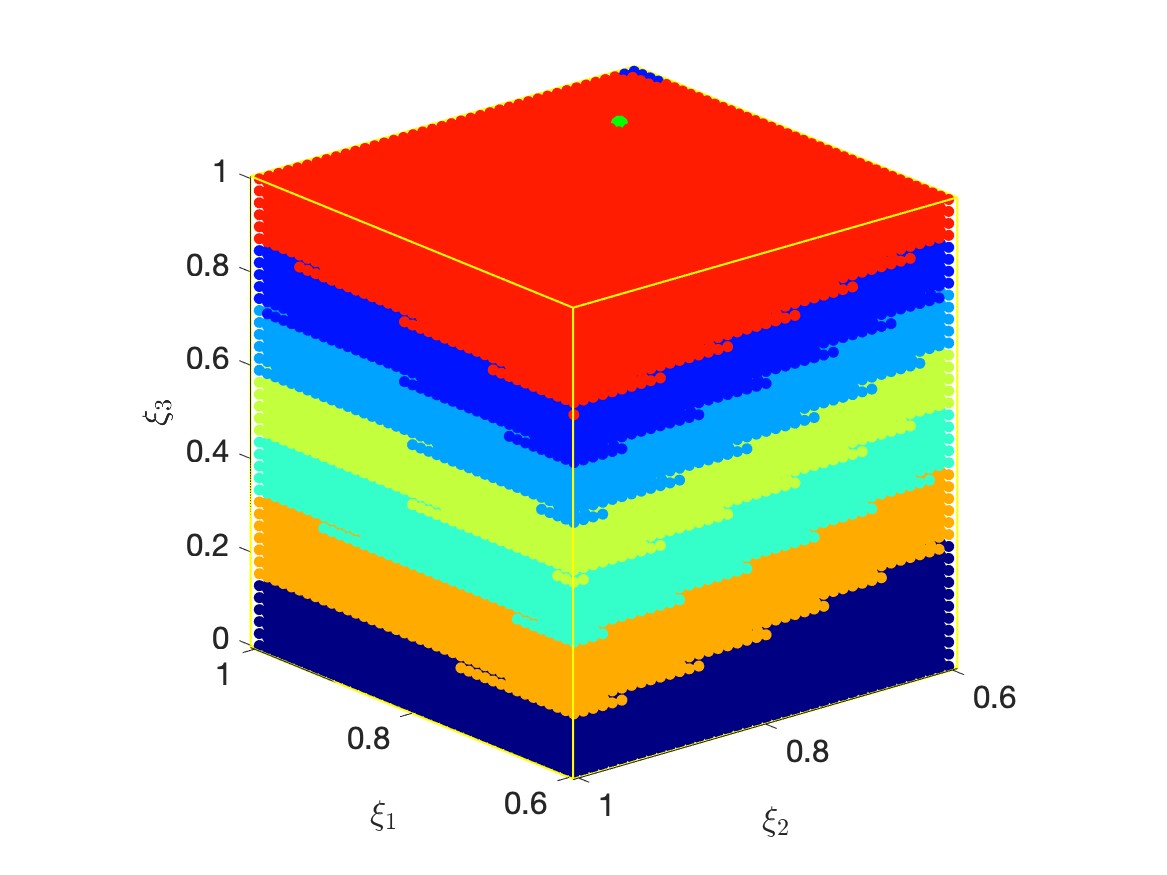}   
    }
    \caption{Visualization of the parameter space: Parameter triples $(\xi_1,\xi_2,\xi_3) = (\xi_\lambda,\xi_A,\xi_\gamma)$ corresponding to dictionary entries in each of the seven subdictionaries are indicated with a different color. The two cubes represent two different view angles, and to facilitate the interpretation, a green fiducial  dot is added on the top of the cube.}\label{fig:cubes}
\end{figure}

Figure~\ref{fig:cubes} shows the parameter space, where each parameter triplet $(\xi_1,\xi_2,\xi_3) = (\xi_\lambda,\xi_A,\xi_\gamma)$ in the sample grid is indicated by a color associated to the corresponding subdictionary. The third parameter seems to be most significantly affecting the shape of the pH curve, as the subdictionaries depend only weakly on the values of the parameters $\xi_1$ and $\xi_2$.

Having the partitioning of the dictionary into the subdictionaries, a criterion is needed for how many feature vectors are included in the code books $\mW^{(j)}$. To find the effective dimensionality of each subdictionary, we consider the first few singular values of each matrix $\mD^{(j)}$. It turns out that the singular values are decaying rapidly, and in all subdictionaries, the ratio of the third and first singular values is below the threshold of $0.001$. We therefore choose $k=3$ for each subdictionary, and compute the NMF decomposition $\mD^{(j)} \approx \mW^{(j)}\mH^{(j)}$ using the IAS-based alternating non-negative least squares solver, see \cite{Waniorek} for details.

To estimate the dictionary compression error, we fix the sample size of the compression errors to $K_i = 7\,000$ for each subdictionary.
For each subdictionary, we draw from the regular grid in the parameter space a random parameter vector $\widehat\xi$ that is associated to the $i$th dictionary. To generate a data vector that presumably represents the $i$th dictionary, we perturb the components of $\widehat \xi$,
\[
 \xi_j = \widehat\xi_j + s \Delta\xi_j, \quad s\sim{\rm Uniform}([-1/2,1/2]), \quad j=1,2,3,
\]
where $\Delta\xi_j$ is the discretization interval of the $j$th component in the parameter space. Using the forward map $F$ from the parameter space to the data base, we compute a feasible data vector, $d = F(\xi)$, and find its corresponding low rank approximation in terms of the feature vectors of the $i$th dictionary, to finally get the compression error sample (\ref{error sample}). We repeat the process $K_i$ times, and obtain estimates for the mean and covariance of the compression error.

We point out here that the generation of the dictionaries and their pre-processing may be time-consuming, and it should be seen as the equivalent of the learning process of neural networks.

The algorithm is tested on the following three experiments.

\begin{enumerate}
\item In the first experiment, we generated the data 
using the generative model in the article \cite{Filtro} 
\begin{equation}
(\xi_\lambda,\xi_A,\xi_\gamma) = (0.9,0.8,0.8) \, ,
\end{equation}
which translates into physical parameter values
\begin{equation}
(\lambda,A^0,\gamma)
=(30.78\, \mu{\rm m}/{\rm s},
16,10^{-4.2}\, \mu{\rm m}/{\rm s}).
\end{equation}
\item In the second experiment, the membrane permeability is tuned down while both the carbonic anhydrase action and the quench factor are tuned up,
\begin{equation}
(\xi_\lambda,\xi_A,\xi_\gamma)
=(0.8,1,2/3), \quad\mbox{ or }\quad
(\lambda,A^0,\gamma) = 
			(27.36\, \mu{\rm m}/{\rm s},
			20,
			10^{-4}\, \mu{\rm m}/{\rm s})
\, .
\end{equation}
\item In the third experiment, the generative model parameters are
\begin{equation}
(\xi_\lambda,\xi_A,\xi_\gamma) = 
(0.8,0.8,0.7667), \quad\mbox{ or }\quad
(\lambda,A^0,\gamma)
=(27.36\, \mu{\rm m}/{\rm s},16,
10^{-4.15}\, \mu{\rm m}/{\rm s})\, .
\end{equation}
\end{enumerate}
To simulate the limited measurement precision, the data was rounded off with the precision of 0.02 pH units.
To identify the best explaining subdictionary for each of the test cases, we run the IAS algorithm  with the sparsity parameter $\eta = 0.001$, and the stopping criterion for the outer iteration $\| \theta^{k+1} - \theta^{k} \| < \text{tol}_\theta = 0.001$. The maximum number of IAS iterations allowed is 1\,000.

The least squares problem with non-negativity constraints within the IAS algorithm is solved by using the CGLS algorithm with a projection onto the non-negative cone. The stopping criterion for the iterations is a Morozov-type discrepancy condition, i.e., the iterations are stopped when the current $x$ vector satisfies
\[
\left\| 
\left(\mC_{\rm DCE}^{(i)}\right)^{-1/2}\left(\mW^{(i)} x - b - \mu_{\rm DCE}^{(i)})\right) 
\right\|^2 \approx T,
\]
the right hand side being equal to the expected square of whitened noise with $T=1\,001$ independent components. 

The classification step of our procedure finds that in experiments one and three, the best fit for the data is given by  the class $2$, while in the second experiment, the best explaining class is $3$.

Once the best subdictionary is identified, the representative dictionary entry is found by solving the problem (\ref{eq:LS_reconstruction}) using the sparsity-promoting non-negative IAS algorithm, with hyperparameter values $\eta = 0.03$ and ${\rm tol}_\theta = 0.001$, with the maximum number of IAS iterations set to 250. In the CGLS step, we use the standard deviation $\sigma=10^{-5}$ as a target for solving the least squares step.

\begin{figure}
	\centering{\includegraphics[width=\textwidth]{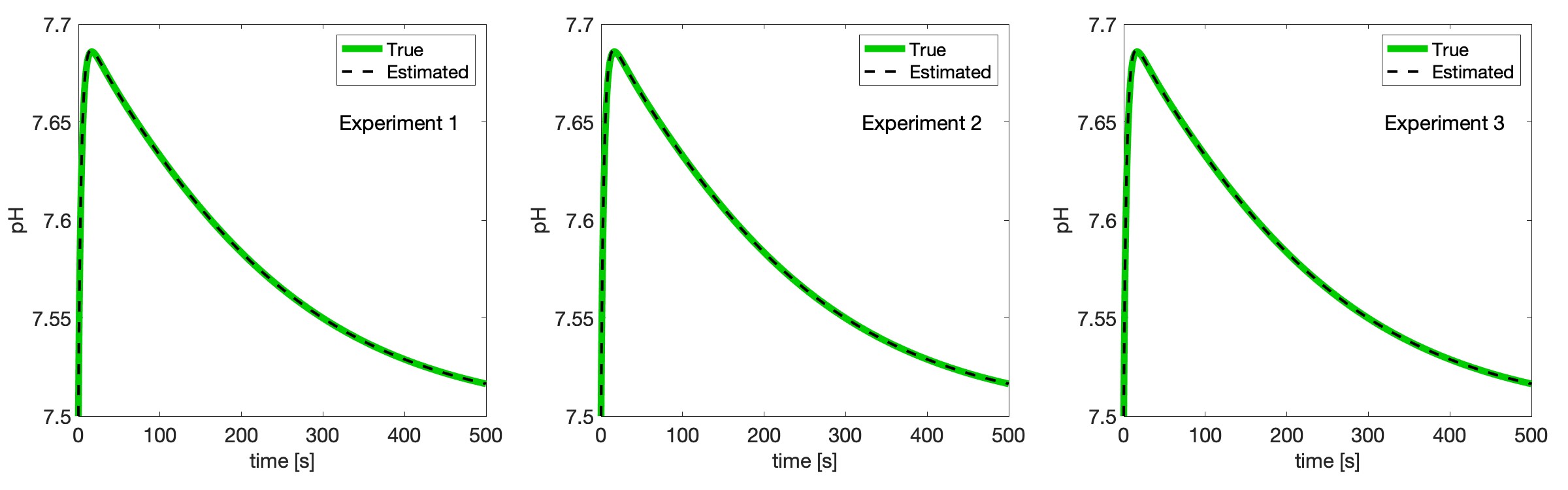}}
	\caption{The pH curves of the three computed experiments. The solid green curve, labeled as "true', represents the computed output with parameter values that were used to generate the synthetic data, and the dashed black curve, labeled as "estimated", is the pH response computed with the parameter values estimated by the dictionary learning algorithm.}
	\label{fig:results}
\end{figure}

The results of the parameter estimation process in the three experiments are shown in Table~\ref{tab:results}, indicating a good performance of the algorithm. To see how well the estimated parameter values reproduce the pH curve, in Figure~\ref{fig:results} the forward model outputs are shown with both the true and the estimated parameter values. We see that in all three cases, the outputs are almost indistinguishable.

\begin{table}
    \centerline{
\begin{tabular}{|c|c|c|c|c|c|c|}
\hline
& \multicolumn{2}{|c|}{Experiment 1}& \multicolumn{2}{|c|}{Experiment 2}& \multicolumn{2}{|c|}{Experiment 3} \\
& True & Est  & True & Est& True & Est\\
\hline
$\lambda$ $[\mu{\rm m}/{\rm s}]$ &30.78 & 30.84& 27.36 & 27.43 & 27.36 & 27.29\\
$A^0$ &16 & 16.29 & 20.0 & 19.36 & 16.0& 17.31 \\
$\gamma$ $[\mu{\rm m}/{\rm s}]$  &$10^{-4.2}$ & $10^{-4.19}$ & $10^{-4.0}$ & $10^{-4.01}$ & $10^{-4.15}$ & $10^{-4.12}$ \\
\hline
\end{tabular}
}
\caption{Results of the parameter estimation process in the three experiments.}\label{tab:results}
\end{table}

Finally, we compare the computingtimes of the proposed algorithm and the particle filter algorithm \cite{Filtro}. While both algorithms can be run without problems in a personal computer environment, for the timings of the test runs reported below the CINECA Leonardo HPC environment was used. In Table~\ref{tab:timing}, shows the timing for each subtask of the dictionary matching algorithm is listed separately. The full calculation using the dictionary of size 64000 atoms from end-to-end takes less than 2 and a half hours, and with precomputed dictionary, the matching process requires about three minutes. For comparison, the parameter estimation using the particle filtering approach with 4800 particles on the same hardware requires from 4 hours 45 minutes to 5 hours using 112 cores, or from 7 hours 30 minutes to 8 hours using 56 cores. The comparison shows that even if the full dictionary is created from scratch, the computing time is reduced by a factor of two or more. Using a precomputed dictionary, the speedup from hours to few minutes is dramatic.

\begin{table}\label{tab:timing}
\centerline{
\begin{tabular}{l|r| l}
Task & \# cores & time \\
\hline
Dictionary creation & 40 & 1 hour 10 minutes \\
Distance matrix & 112 & 33 minutes \\
$k$-medoid clustering & 112 & $<$ 10 minutes\\
Rank estimation & 112 & seconds \\
NMF factorization & 112 & 4 minutes \\
Error density estimation & 112 & 23 minutes \\
Dictionary matching & 8 & 3 minutes
\end{tabular}
}
\caption{Timing of the subtasks of the dictionary matching algorithm with a dictionary comprising 64000 atoms. The computations were performed on the CINECA Leonardo HPC having 112 cores/node: Each node has two Intel Sapphire Rapids Intel Xeon Platinum 8480+ processors of 56 cores.}
\end{table}

\section{Conclusions}

In this article, we propose a layered dictionary matching algorithm for the estimation of model parameters of a computational model describing the pH response due to gas transport through a cell membrane in a presence of a measurement device that itself has an effect on the microenvironment.
The main characteristics of this model is that while the parameter space is low-dimensional, the forward model is complex and computationally intensive. Moreover, the non-linearity of the forward model makes the analysis of optimization-based parameter estimation algorithms difficult, while, e.g., sampling-based solutions of the inverse problem may be too time consuming. In
\cite{Filtro}, the authors considered the inverse problem in the framework of Bayesian particle filtering techniques which avoids the time-consuming forward solutions of the full model over the experimental time window, requiring only local solutions over short time steps. Compared to the particle filtering method with $4800$ particles and $1001$ measurements at $2 {\rm Hz}$, the computing times of the proposed dictionary learning method is, given similar hardware, approximately from one half to one third, depending on the hardware.
This comparison takes into account the preprocessing time comprising the creation of the dictionary, clustering, extraction of the features and dictionary compression error estimation. We point out that once the preprocessing phase has been completed, solution of the inverse problem is very fast, thus improving the performance by orders of magnitude. While the filtering approach could be made faster, e.g., by resorting to other approximate Bayesian filtering methods, the Ensemble Kalman Fitering (EnKF) or Unscented Kalman Filtering (UKF) in particular being viable candidates, the performance of the dictionary learning technique with precomputed dictionaries remains hard to beat.

Beyond the comparison of the computing times, the direct comparison of the filtering algorithm and the dictionary learning approach is not straightforward because of the notable differences in the approaches. In particular, the particle filtering approach allows
uncertainty quantification, which is built in the Bayesian  filtering paradigm. Moreover, the particle filtering approach allows the model parameter to change over time, while in the dictionary learning method, a static parameter model is necessary. In the present application, such assumption may be natural, while in general, allowing time dependency of the model parameters may add flexibility. Extending dictionary learning techniques for non-static parameter models is an interesting challenge.

The technique proposed in this paper appears useful in the case in which the estimation of parameters must be carried out with limited computational resources, in a relatively short time, or when several estimates must be run in parallel starting with the same dictionary. The obtained estimates appear to be relatively accurate, in fact, the predicted outputs with estimated parameters are indistinguishable from the original data, as Figure~\ref{fig:results} suggests. In the numerical experiments with simulated data, the noise was restricted to a coarse truncation error noise. To have a better understanding of the sensitivity, a careful analysis of the measurement noise and inaccuracies should be carried out. A limitation of the proposed algorithm is its reliance on dense sampling of the parameter space which was low-dimensional in the present example. In high-dimensional parameter models, the parameter space sampling requires further thoughts, e.g., by relying multilevel sampling or  Monte Carlo sampling with some extra prior information that guides the sampling process.

\section*{Acknowledgements}

The work of DC was partly supported by the NSF grants DMS 1951446 and  DMS 2513481, and that of  ES by the NSF grant DMS 2204618 and DMS 2513481.

\end{document}